\documentclass[journal, 9pt]{IEEEtran}
\IEEEoverridecommandlockouts
\usepackage{cite}
\usepackage{amsmath,amssymb,amsfonts}
\usepackage{algorithmic}
\usepackage{graphicx}
\usepackage{textcomp}
\usepackage{xcolor}
\usepackage{caption}
\usepackage{enumerate}
\usepackage{subcaption}
\usepackage[compact]{titlesec}
\usepackage{multirow}
\usepackage[ruled,vlined]{algorithm2e}
\begin{document}

% \titlespacing{\section}{0pt}{0pt}{0pt}
\titlespacing{\section}{5pt}{5pt}{5pt}
\setlength{\abovedisplayskip}{5pt}
\setlength{\belowdisplayskip}{5pt}

% \title{The Value of Operation Coordination to EV Fleet Aggregators}
\title{The Value of Operational Coordination for EV Fleet Aggregators }
\author{\IEEEauthorblockN{Oluwaseun Oladimeji, David Pozo,}
\IEEEauthorblockA{\textit{Skolkovo Institute of Science and Technology}\\ }
\and
\IEEEauthorblockN{Yury Dvorkin,}
\IEEEauthorblockA{\textit{New York University}\\}
\thanks{
\hrule \vspace{.1cm}
This work was partially funded by the Skolkovo Institute of Science and Technology as a part of the Skoltech NGP Program (Skoltech-MIT joint project).}
\vspace{-0.6cm}
}
\maketitle
\IEEEpubidadjcol

\begin{abstract}
% Integration of energy systems such as power and gas networks as well as power and thermal networks can bring significant advantages to system security, reliability, and decreased emissions. Another alternative for sector coupling with large potential benefits is the power and transportation network. This is largely due to the increasing electric vehicle (EV) usage and their demand on the electric power grid. Additionally, costs of production and operation of EVs and battery technologies are consistently decreasing while tax credits on EV purchase and usage are being offered to users in developed countries. The power grid is also undergoing major updates and changes with focus on ensuring a greener environment. These factors influence our paper.
% We present a new operation model for an EV-grid integrated system which incorporates a set of aggregators (who own fleet of EVs) with part access to the distribution grid. Cooperative game theory is then used in modelling the system behavior. The Core is used in describing the stability of the interaction between these aggregators and Shapley value is used in allocating costs to them. Results obtained show the benefit of cooperation which could either lead to overall reduced energy consumption, reduced operation cost of EVs and the distribution grid, and in some cases, the available extra cost for strengthening the transport and electric grid infrastructures.
The integration of energy systems such as electricity and gas grids and power and thermal grids can bring significant benefits in terms of system security, reliability, and reduced emissions. Another alternative coupling of sectors with large potential benefits is the power and transportation networks. This is primarily due to the increasing use of electric vehicles (EV) and their demand on the power grid. Besides, the production and operating costs of EVs and battery technologies are steadily decreasing, while tax credits for EV purchase and  usage are being offered to users in developed countries. The power grid is also undergoing major upgrades and changes with the aim of ensuring environmentally sustainable grids. These factors influence our work.
We present a new operating model for an integrated EV-grid system that incorporates a set of aggregators (owning a fleet of EVs) with partial access to the distribution grid. Then, the Cooperative Game Theory is used to model the behavior of the system. The Core is used to describe the stability of the interaction between these aggregators, and the Shapley value is used to assign costs to them. The results obtained show the benefit of cooperation, which could lead to an overall reduction in energy consumption, reduced operating costs for electric vehicles and the distribution grid, and, in some cases, the additional monetary budget available to reinforce the transmission and grid infrastructures.
\end{abstract}

\begin{IEEEkeywords}
Aggregators, Cooperative Game Theory, Electromobility, EV, Aumann Game
\end{IEEEkeywords}

\section{Introduction}
The numerous advantages of EVs have made them a good choice for grid power management, voltage regulation and power quality improvement \cite{das2019electric}. The possibility of operating EVs as energy storage devices, components of virtual power plants and as scheduled loads for peak shaving has made them increasingly popular in power system studies. Coupled with these advantages is their replacing the internal combustion engine vehicles for climate change mitigation.

However, some negative effects have also been associated with EV integration with the power grid. EV loads are nonlinear because human behavior is unpredictable, thus, charging behavior is also largely unpredictable. Moreover, increasing usage of fast chargers means large amount of power can be drawn in a short time which further leads to the challenge of stability of the power grid \cite{gomez2003impact}. Additionally, a large number of EVs leads to a high amount of power requirement from the grid. The enormous amount of real power consumption leads directly to large power loss in the distribution system \cite{fernandez2010assessment}. A coordinated charging scheme can however be a suitable solution to this challenge. As stated earlier, however, human actions are unpredictable and coordination can only be possible on an aggregated scale.

Coordination schemes have been proposed to address some of the challenges stated. Some of the coordinating objectives have included optimization of total energy cost by minimizing the energy losses due to charging/discharging schedules of a large number of EVs and minimizing operation cost \cite{wu2018dynamic}. Taking the distribution network into account, the optimization of EV charging  while considering the distribution network constraints have also been reported \cite{de2014optimal}.

One of the EV charging coordination schemes that has been proposed is valley filling. Offline and online valley filling algorithms with time-dependent optimal power flow problem were used in \cite{chenoptimal,chenvalley}. 
% The authors recognised the challenge posed by large EV fleet on the distribution network if they are stationed in residential areas. The EV charging problem was modelled using a modified version of the optimal power flow problem with a valley-filling objective. 
In \cite{callawaydecen}, a mode of decentralized coordination of charging large populations of EVs for electricity costs minimization was proposed. \cite{karfopoulos2012multi} presented a mode of efficient grid operation while fulfilling the energy requirements of a large number of EVs. In summary, these methods seek to create a profitable scenario for both EV users and the power grid while minimizing power loss and maximizing load factor.

In this paper, we extend our definition of coordination to transcend numerous single-user EVs and we consider another approach similar to integration of power systems. Systems integration has advantages of decrease of global operation cost, systems reliability and environmental benefits.
% and possible increase of the social welfare function.
Integration of electricity and thermal systems in microgrids was proposed in \cite{liu2018hybrid} with minimization of overall operation cost as the objective. An energy management operation model combining bio-gas, solar and wind renewables was presented in \cite{xu2018distributed} while \cite{liu2017distributed} proposed a robust energy management scheme for multiple interconnected microgrids in the real time energy market. While integration of multiple systems might or might not include trade-offs, it is customary to evaluate costs and/or benefits and how these are distributed among the different entities involved in the integration.

Non-cooperative and cooperative game theory frameworks have been proposed in the subject of integration, respectively, with respect to maximizing utility and sharing costs among the different players. Non-cooperative game theory was utilized in \cite{alipour2018multi} to solve an optimization problem involving multiple CHPs with uncertainties in market prices and customer demand. An energy management scheme for coordinated electricity and heat demand response scheme was presented in \cite{liu2017multiparty} using non-cooperative game theory as well.

For cooperative game theory, Shapley value is used in the distribution of the worth of each coalition, either in terms of costs or benefits. In \cite{chics2017coalitional}, cooperative game theory model for several households sharing renewable and energy storage space was presented with the objective of reducing energy consumption costs. A coalitional game based approach for power exchange by different microgrids on a network was proposed in \cite{mei2019coalitional}. An inherent challenge, however, with cooperative game theory is intractability which might arise as the number of entities involved increases.

An important consideration that has not been highlighted in the literature is interoperability of charging stations owned by more than one aggregator. In fact, most studies that have included aggregator have only modelled a single fleet owner with many EVs. Our major contribution, thus, is the formulation of a coalitional-game based model involving multiple aggregators in a distribution grid-transport network integrated system. These aggregators can use each other's charging stations and the value of these coalitions will be important for fleet owners in making decisions whether to cooperate with one another. This coalition value, which is a criteria for decision making by the EV aggregators, is one of the objectives of our investigation. The impact of these coalitions on overall energy consumption and $CO_2$ emissions will also be addressed. We address the problem of intractability by the use of aggregators. The maximum number of fleet owners is not greater than six (6) in all our scenarios, thus, we can enumerate all possible combinations or coalitions and see the interactions between multiple players.

The rest of the paper is organized as follows: The EV to grid integration physical framework will be described in Section \ref{physical system}. The characteristics of the aggregators and DSOs will be highlighted. Section \ref{math framework} describes the cooperative game-based model, routing model for each EV, energy management model for each EV and distribution grid formulation. Results, discussions, and conclusions will be enumerated in Sections~\ref{results} and \ref{conclusion}.

\section{EV-Grid Integration System} \label{physical system}
\begin{figure}
     \centering
     \includegraphics[scale = 0.35]{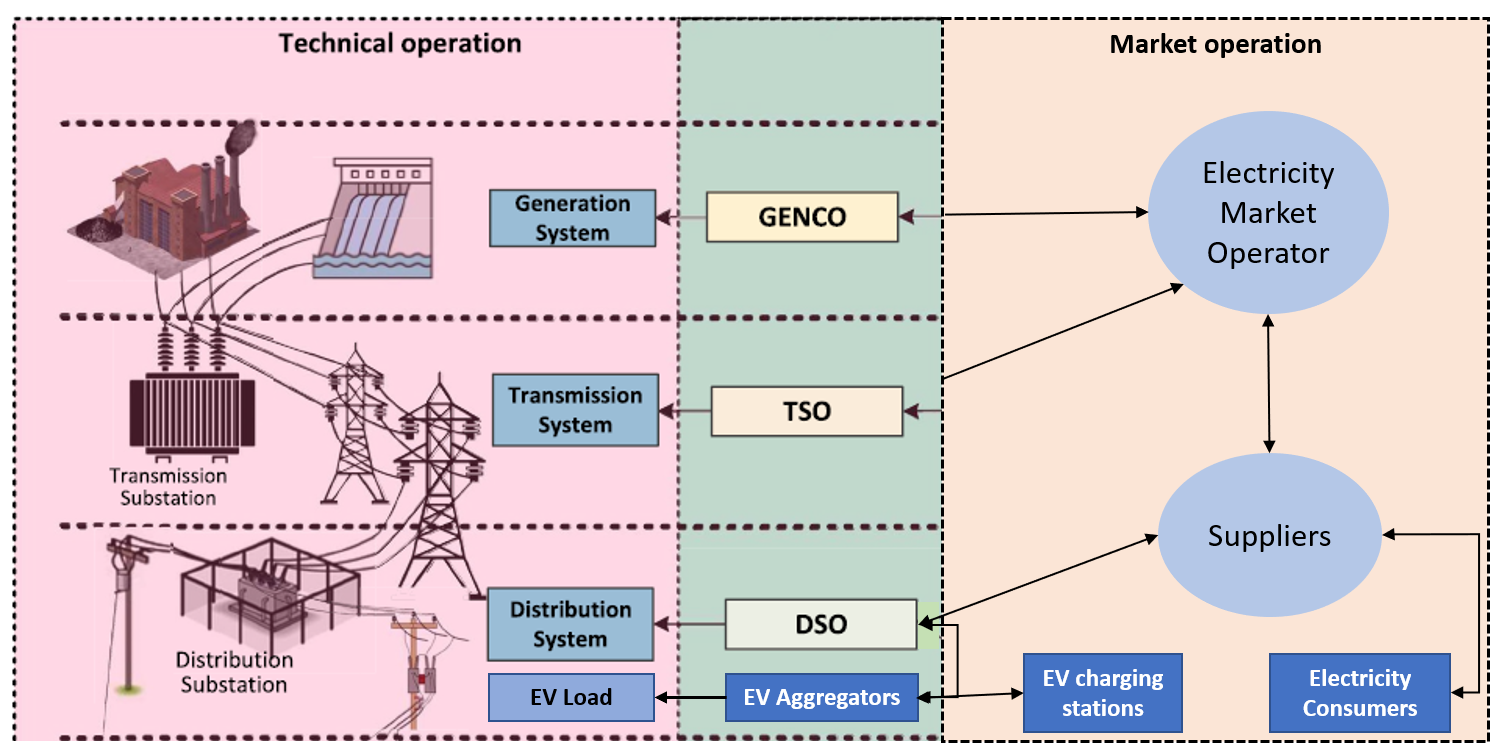}
     \caption{Grid and EV Fleet Coordination in Power System}
     \label{fig:physcial system}
     \vspace{-0.5cm}
\end{figure}

The physical system we propose, adapted from the EV integration framework in \cite{das2019electric}, is shown in Fig~\ref{fig:physcial system}. For our model, aggregators work directly with the distribution system operators, DSO, instead of unilaterally getting directives from DSO or supplying information alone. Other characteristics of the aggregators are:
\begin{itemize}
    \item Control of EV charging stations
    \item Forecast of EV behaviour and transport profiles. This is possible as aggregators function in the role of transport service providers
    \item Access to the localized power consumption data just like the DSO
    \item In the market side, they can act as wholesale agents in the power market
    \item They can also provide ancillary services such as demand response, frequency regulation or act as virtual power plants 
\end{itemize}

These characteristics of the aggregators allow for making autonomous decisions. Likewise, access to the local power consumption data ensures that each aggregator can make decentralized decision of its charging schedule, thus, guaranteeing grid stability. In the mathematical framework described in Section \ref{math framework}, the reader should keep in mind that we have a number of aggregators who are in charge of $N$ number of EVs which are primarily used as ride-sharing cars. These vehicles are not static for long periods of time, thus, there is the possibility of one aggregator requiring the services of charging stations of another aggregator.

\section{Mathematical Framework} \label{math framework}
The overall aim of the cooperative game which subsumes the optimization framework is to efficiently allocate cost of the most efficient coalition to aggregators. This translates to minimizing travel time, energy usage of the EVs and cost of the distribution network operation by efficiently scheduling charging periods.

However, in order to have a benchmark to compare results from our cooperative cases, we present an incomplete information scenario borrowing the idea of the \emph{Aumann model of incomplete information}. The choice of this framework is motivated by some real life scenarios which may not be modelled using extensive form and strategic-form games and incomplete information scenarios fall largely under this category.

\subsection{Model of Incomplete Information}
Modelling our case study as a variant of non-cooperative game is a classic example of games with incomplete information. The situations here are such that aggregators do not have complete information on the environment they face. As a result of the interactive nature of the game, formulating such situations involves not only the knowledge and beliefs of the players, but also the whole hierarchy of knowledge of each player \cite{maschler2013game}. 

This rather infinite knowledge base might be infeasible in evaluating real life scenarios. Studies have shown, however, that most people rarely engage in more than second order common knowledge (about all other players' knowledge) which makes evaluation of scenarios feasible \cite{camerer2011behavioral}. In our case, there are some common knowledge which primarily include the transport and electric grid network.
Here, we utilize the Aumann model of incomplete information with beliefs for our formulation. The solution here will arise from players settling on just the prior information received before facing the game. A formal definition of the Aumann model consists of $$(N, Y, (F_i)_{i \in N}, s, \mathbf{P})$$ where:
\begin{itemize}
    \item $N$ is a finite set of players;
    \item $Y$ is a finite set of states of the world;
    \item $F_i$ is a partition of $Y$, for each $i \in N$;
    \item $s: Y \longrightarrow S$ is a function associating a state of nature to every state of the world;
    \item $\mathbf{P}$ is a probability distribution over $Y$ such that $P(\omega) > 0$ for each $\omega \in Y$; $\omega$ is a state of the world.
\end{itemize}

$N$ represents EV aggregators in our case. The state of nature is the parameter with respect to which there is incomplete information: for an e-mobility example, the route through which EV $k$ owned by aggregator I takes passengers from origin $i$ to destination $j$. Thus, $F_I = \{x_{ijk} \, \, \forall i,j,k \}$: all routes from origin to destination by vehicles owned by aggregator I. The same can be defined for other aggregators as required.

% The state of the world, $\omega$, includes, in addition, the common knowledge structure of the players regarding the state of nature; for example, the state of the world corresponding to the above state of nature goes further to say aggregators II and III know this route taken by vehicle $k$ to visit destination $j$ from origin $i$. The state of the world is thus dependent on the characteristics of the belief spaces (i.e. common prior). The assumption of common prior being shared by all aggregators is however a strong assumption and in many games, its existence is questionable. This then leads to the inclusion of probability distribution $\mathbf{P}$ which gives the probability that any specific state of the world in $Y$ is the true one.

To enumerate this model, we present two games with different belief spaces; the beliefs we employed are chosen based on current system practices. For our particular subject, consistent beliefs held by all aggregators are (i) knowledge of transport and electric grid (i.e. all EV aggregators in New York City have perfect knowledge of the transport network and the power grid of New York City) and (ii) knowledge of relative locations of other aggregators (i.e. aggregator A knows the location of aggregator B's assets, depot etc.)
The difference between our games occurs in cases when (a) users are \textit{aggregator-aware}, i.e., much like the current transportation service providers landscape, users can decide which aggregator(s) to call using their mobile applications or other available access points; (b) users are \textit{aggregator-blind} which means users do not choose which aggregator provides their transportation requirements. Variants of these two games were outlined in \cite{fisk1984game}.

\subsection{Cooperative Games Allocation Rules}
% The question of inter-fleet coordination to use each other's resources for service provision has been raised in the beginning of this section. 
We propose modelling a coalition game where the Core is analysed to determine if there are opportunities for cooperation and then costs or benefits are shared to each player by Shapley value computation. The Shapley value is a solution concept used for fair and efficient distribution of gains or costs to several players working in a coalition \cite{shapley1953value}. For this work, we evaluate the cost that each aggregator has to pay and thus, we will formulate a cost sharing game.

\subsection{The Core}
Let $(N;c)$ be the cost-sharing coalitional game where $N$ is the number of EV aggregators, and $c$ is the coalitional function.

The Core is used in describing the solution set of all stable cost allocations obtained. For the problem(s) presented, there are two obvious possibilities which are defined by stability and fairness. Either all players will form the grand coalition where they all cooperate or they deviate and form sub-coalitions. The scenario that an adjudicator wants is for every player to be a part of the grand coalition and reduce the possibility of deviation from this. The grand coalition translates to the least operation cost for EV and grid in our case. Given $x$ as the allocation vector cost for all players, the Core is formulated as
\begin{subequations} 
\begin{IEEEeqnarray}{r'll}
    X(N;c) = \Big\{x \in \mathbb{R}^N: 
    & x_i\leq c(i), &  \forall{i \in N} \label{eq. ind_rationality} \\
    & \sum\limits_{i \in N} x_i = c(N),  \label{eq: efficiency}  \\
    & \sum\limits_{i \in S} x_i \leq c(S), & \: \forall{S \subseteq N} \Big\}.  \label{eq: rational subcoalition}
\end{IEEEeqnarray}
\end{subequations}

The Core is termed the result of rational coalition \cite{maschler2013game}, meaning that: 
(i) no player obtains a higher cost than what he could get on his own, also known as \textit{individual rationality principle} \eqref{eq. ind_rationality};
(ii) the sum of the costs paid by each individual player is equal to the sum of the cost of the grand coalition, also known as \textit{game efficiency principle} \eqref{eq: efficiency}; and 
(iii) for every sub-coalition, $S$, formed which is not the grand coalition, no player $i$ should pay more than it would if acting alone, also known as \textit{coalitionally rational principle} \eqref{eq: rational subcoalition}. 
% The sum of each individual's costs should not be less than cost of grand coalition so that cooperation efforts will not break down and costs should not be greater such that there would be no wastage. The Core is defined essentially to fulfill the principles of:
% (1) efficiency
% \begin{IEEEeqnarray}{l}
%      \IEEEeqnarraynumspace \IEEEyesnumber \label{eq: rational coalition}
% \end{IEEEeqnarray}
% and (2) coalitional rationality
% \begin{IEEEeqnarray}{l}
     
% \end{IEEEeqnarray}
% i.e.  This is a rational argument and it gives players more incentive to cooperate since they expect to incur no greater cost when they cooperate. Also, the cost of the grand coalition $c(N)$ is lower than every sub-coalition $c(S)$ formed.

We should note that it is not every instance where there are workable agreements between players and games can fail at those instances. However, when the Core is non-empty, it aids in determining the stability of the coalition(s) while the Shapley value, a single-valued solution concept, is used for allocating costs to each player who participates in the game.

\subsection{Shapley Value}
The Shapley value, which we use for cost sharing among players in here, is computed by 

\begin{IEEEeqnarray}{ll}
    Sh_i(N;c) = & \sum\limits_{S \subseteq N \setminus \{i\}} \frac{|S|!  (|N|-|S| -1)!}{|N|!} \Big(c\left( S\cup\{i\} \right) - c\left(S\right)\Big) \nonumber 
    \vspace{-0.2cm} \\&  \IEEEyesnumber  \label{eq: shapley value} 
\end{IEEEeqnarray}
\noindent where $|N|$ is the total number of players in the game and the sum ranges over the subsets $S$ of the coalitions which does not contain the player $i$.
%The sum ranges over all $|S|!$ orders $R$ of the players where $S$ is the set of players in $N$ which precede $i$ in the order $R$. $c(S)$ represents the value function for all subsets of the game. E.g., for a two player game, there are only two instances; with opportunity for either player to play first. The orders are \{1,2\} and \{2,1\} while the game subsets include \{1\}, \{2\} and \{1,2\}.
The term $c\left( S\cup\{i\} \right) - c\left(S\right)$ represents the marginal contribution of player $i$ to the coalition $S$. This value is the solution of the EV transport routing problem described below.

\subsection{EV Transport Routing}
The transport routing problem answers the question of optimal assignment of routes to a fleet of vehicles based on the needs of customers in different locations in a transportation network. The main objective of transport routing is the minimization of the total distance traversed by all vehicles, although other alternative objectives are possible.

For EVs, an efficient battery management system needs to be modelled to account for energy usage. The objective here is minimization of the overall cost of energy used and this includes the reduction of time spent in charging stations outside the depot. Additionally, the distribution grid is coupled in our formulation because our approach leverages the benefits of the co-optimization of the electric-transport integrated system. This co-optimization ensures the stability of voltage variations in distribution grid is within acceptable bounds and satisfaction of network constraints. A linear representation of the distribution grid, LinDistFlow \cite{baran1989optimal}, is utilized here. This also allows easy coupling of the energy drawn by the EV to the local consumption at the charging station. Full formulations for EV routing, energy management and distribution grid are outlined in online appendix \cite{onlinecomp}.

\subsection{Overall Framework}
\begin{figure}[ht]
     \centering
     \includegraphics[scale = 0.7]{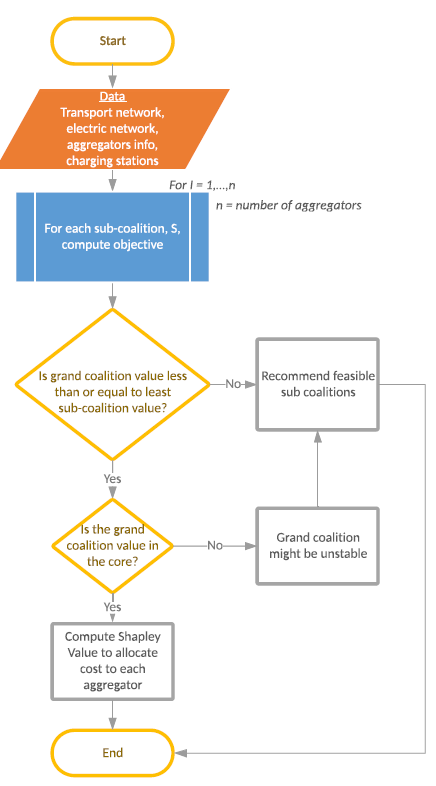}
     \caption{Flow chart}
     \label{fig:flowchart}
\end{figure}

Finally, our goal is to evaluate the allocation of operation cost of an EV fleet and distribution grid when the aggregators cooperate using the flowchart in Figure~\ref{fig:flowchart}. This is achieved by enumerating possible coalitions by aggregators and then computing the share that accrues to them in a successful coalition. The successful coalition is thereafter evaluated for the assignment of routes and charging schedules.

\section{Results and Discussions} \label{results}

\subsection{Transportation Network Data}
The case study selected for validating our proposed methodology is based on the Downtown Manhattan transportation network. Drop-off and pick-up data is taken from taxi trip records by the New York City Taxi and Limousine Commission for 2018 \cite{nyctaxidata}. The data set includes pick-up and drop-off dates, times, locations, trip distances and driver-reported passenger counts.
% Table~\ref{table:data snippet} shows a snippet from this New York City data set where PULID and DOLID are respectively the pick-up location ID and drop-off Location IDs, which indicate exact latitudes and longitudes of these locations.
%{Historical data (partly represented in Table~\ref{table:data snippet}) was analysed and it showed that drop-offs and pickups follow a Gaussian distribution.}

% \begin{table}[ht]
% \begin{center}
% \caption{Representative Data Snippet}
% \begin{tabular}{c c r r r r} 
%  \hline
%   Pickup & Dropoff & \multicolumn{1}{c}{Trip} & \multicolumn{1}{c}{PULID} &  \multicolumn{1}{c}{DOLID} & \multicolumn{1}{c}{Passenger} \\ 
%  date-time & date-time & \multicolumn{1}{c}{distance} &  &  & \multicolumn{1}{c}{count} \\ [0.5ex]
%  \hline\hline
%     12/31/2018  & 1/1/2019  & 16.6 & 162 & 26 & 1\\     
%     11:58:45 & 12:37:04 & & & & \\
%     \hline
%     12/31/2018  & 1/1/2019 	& 3 & 144 & 170 & 2 \\     
%     11:58:44 & 12:16:54 & & & & \\
%     \hline
%     12/31/2018 	& 1/1/2019 	& 1.8 & 158	& 68 & 1 \\ 
%     11:58:43 & 12:13:49 & & & &\\
%     \hline
%     12/31/2018 	& 1/1/2019  & 1.1 & 141	& 237 & 1 \\
%     11:58:38 & 12:03:29 & & & & \\
%     \hline
%     12/31/2018 	& 1/1/2019 	& 6.4 & 162 & 13 & 1 \\
%     11:58:37 & 12:14:04 & & & & \\
%     [1ex]
%  \hline
% \end{tabular}
% \label{table:data snippet}
% \end{center}
% \end{table}
Using latitudes and longitudes, we divided the data into distinct locations such as commercial, industrial, and residential regions. With this division, we were then able to track movement patterns at different times of the day. For example, people generally move from residential areas to commercial and industrial centers in the mornings and vice versa in the evenings.

\begin{figure}[ht]
     \centering
     \begin{subfigure}[t]{0.425\columnwidth}
         \centering
         \includegraphics[width=\columnwidth]{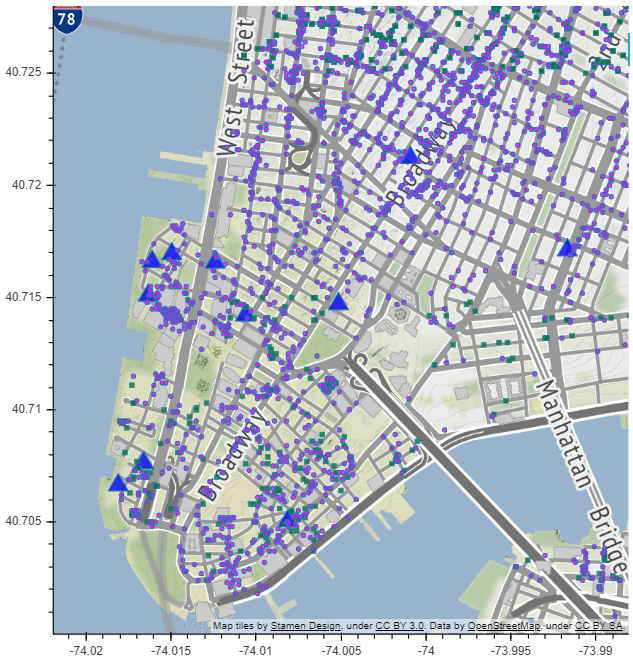}
         \caption{Passenger concentration and EV charging stations in Downtown Manhattan between 6am-10am, 13th August, 2018}
         \label{fig:Downtown manhattan}
     \end{subfigure}
     \hfill
     \begin{subfigure}[t]{0.425\columnwidth}
         \centering
         \includegraphics[width=\columnwidth]{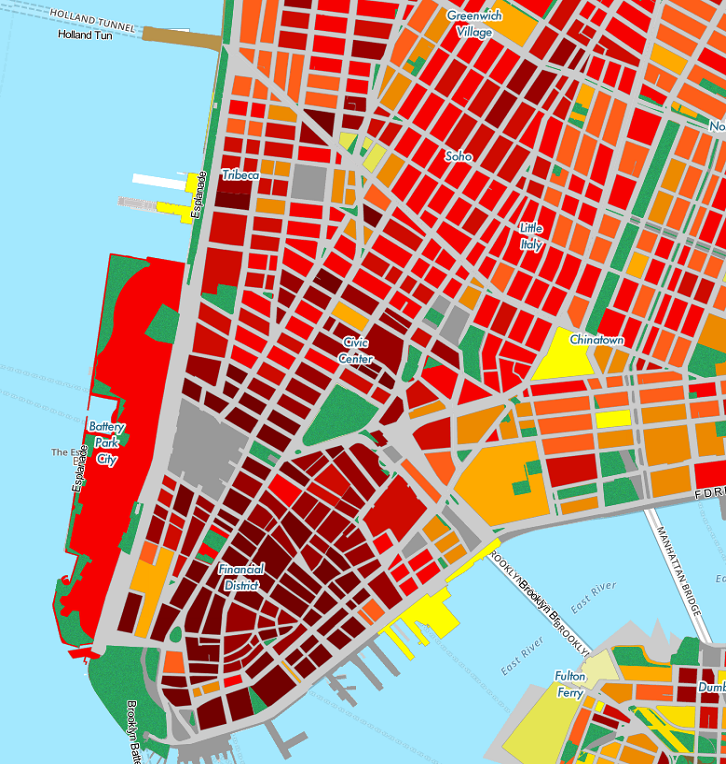}
         \caption{Downtown Manhattan energy consumption by area}
         \label{fig:Downtown manhattan power}
     \end{subfigure}
        \caption{Downtown Manhattan transport and energy demands}
        \label{fig:Downtown manhattan Transport and Energy Consumption}
\end{figure}

Fig.~\ref{fig:Downtown manhattan} demonstrates a geographical snapshot of the demand data, which is at the Downtown Manhattan grid, where green squares represent originating nodes, red circles represent destination nodes and large blue triangles are EV charging stations.
% \begin{figure}
%   \centering
%   \includegraphics[scale=0.5,valign=t]{example-image-a}
%   \includegraphics[scale=0.3,valign=t]{example-image-b}
%   \caption{A caption\label{fig:scaled_diss}}
% \end{figure}

%%%%%%%%%%%%%%%%%%%%%%%%%%%%%%%%%%%%%
\subsection{Electric Network Data}
For local electricity consumption, we obtained information from \cite{howard2012spatial} and New York City Mayor's Office of the Long-Term Planning and Sustainability \cite{nycSIRR}.
Electricity consumption at each block and lot level was obtained from \cite{howard2012spatial} and consumption patterns are shown by the different shades of red in Fig.~\ref{fig:Downtown manhattan power}. Darker red colors show regions of high power consumption while lighter/yellow colors show areas of low energy usage. Voltage ratings of each charging station location and reactance between stations were estimated similarly to \cite{acharya2020public}.

% For local electricity consumption in each area where charging stations are located, we obtained information from \cite{howard2012spatial} and New York City Mayor's Office of the Long-Term Planning and Sustainability \cite{nycSIRR}.
% Electricity consumption at each block and lot level was obtained from \cite{howard2012spatial} and consumption patterns are shown by the different shades of red in Fig.~\ref{fig:Downtown manhattan power}. Darker red colors show region of high power consumption while lighter/yellow colors show areas of low energy usage.

% In Table~\ref{table:electrical data}, we show local electricity consumption of some areas in Downtown Manhattan. Here voltage ratings of each charging station location and reactance between stations were estimated similarly to \cite{acharya2019public}. 

% \begin{table}[ht]
% \begin{center}
% \caption{Distribution Grid Data}
% \begin{tabular}{l r r} 
%  \hline
%   \multicolumn{1}{c}{Block ID} & \multicolumn{1}{c}{Land Area $[m^2]$} & \multicolumn{1}{c}{Electricity Use [kWh]} \\ [0.5ex] 
%  \hline\hline
%     16      & 370007    & 322379    \\     
%     70      & 6666 	    & 30545     \\     
%     140	    & 6541 	    & 4476      \\ 
%     171	    & 7057      & 5622      \\
%     175     & 9479 	    & 8081      \\ [1ex]
%  \hline
% \end{tabular}
% \label{table:electrical data}
% \end{center}
% \end{table}
%%%%%%%%%%%%%%%%%%%%%%%%%%%%%%%%%%%%%
%               Result              %
%%%%%%%%%%%%%%%%%%%%%%%%%%%%%%%%%%%%%

We first present Case A comprising 2 scenarios which serves as our benchmark before we describe the cooperative games comprising 3 to 6 aggregators.
\subsection{Case A: Model of Incomplete Information}
The Aumann game is defined as follows:
\begin{itemize}
    \item The set of aggregators $N$ = \{Agg. I, Agg. II, Agg. III \}
    \item The set of states of the world $Y$ = \{$\mid$routes to satisfy all passengers$\mid$\}
    \item The information partitions of the players are network dependent and this is based on which part of the grid customer requests originate.
    \item The common prior we select as $\frac{1}{3}$ for each state of the world. i.e., there is a 1 in 3 chance that any of the aggregators will be chosen to transport a passenger from an origin to a destination.
\end{itemize}
% The overall objective of minimization of transport and energy cost remains the same.

\subsubsection{Case A1. Aggregator-aware Transport Service}
The first case is for \emph{aggregator-aware users}, that decide by themselves about the aggregator that provides the service. Deduction from real life cases reveals that users choose which transport service provider they want to use at the particular instance. Either the public transit or one of the numerous private ride sharing providers.
The energy usage and cost incurred by each aggregator (obtained by solving the Case A model with model from \cite{onlinecomp}) is shown in Table~\ref{table:3 aggs caseA1}. The values in the table will be compared with the results obtained from cooperative game framework at the end of this section.

\begin{table}[ht]
\begin{center}
\caption{Case A1 Results Summary}
\begin{tabular}{r || c | r}
    \hline
     Aggregator & Value Fnc(\$) & E (kWh)   \\ 
    \hline\hline
        \{1\}     & 2.137 & 8.633  \\    \relax   
        \{2\}     & 4.323 & 17.466 \\    \relax  
        \{3\}     & 2.373 & 9.588  \\   
        [1ex]
    \hline
\end{tabular}
\label{table:3 aggs caseA1}
\end{center}
\vspace{-0.5cm}
\end{table}
As noted in the model, each aggregator has a 1 in 3 chance to be chosen for the routing. Thus, least cost route amongst aggregators cannot be ascertained in this case. This choice, in effect, could lead to choosing the farthest provider at some point. However, this choice of farthest service provider is as likely as the choice of the closest provider. Further, interactions between aggregators could, over time, lead to an equilibrium of efficient operation by the aggregators.

\subsubsection{Case A2. Aggregator-blind Transport Service}
This is why we decided to carry out a similar analysis but with \emph{aggregator-blind users}. It means users/passengers cannot decide which aggregator provides their service. This could be restrictive but, owing to the extent of unknowns that might be encountered, could lead to some interesting inferences. 
As seen in Table~\ref{table:3 aggs caseA2}, aggregator 2 participated only once in the routing while all other routes were made by aggregators 2 and 3. This makes sense since all aggregators now have access to all demands at any time and each could decide whether offering a particular service would be optimal or otherwise.
\begin{table}[ht]
\begin{center}
\caption{Case A2 Results Summary}
\begin{tabular}{r || c | r}
    \hline
     Aggregator & Value Fnc(\$) & E (kWh)   \\ 
    \hline\hline
        \{1\}     & 5.281 & 21.337  \\    \relax   
        \{2\}     & 0.719 & 2.906   \\    \relax  
        \{3\}     & 2.833 & 11.444  \\   
        [1ex]
    \hline
\end{tabular}
\label{table:3 aggs caseA2}
\end{center}
\vspace{-0.5cm}
\end{table}
This is substantially different from the case A1 where aggregator 2 got about half the total value function as well as expended half of the energy consumed. While this case 2 might not be feasible due to recent practises, it is an option which might give results close to coalition efforts as described below.

\subsection{Case B: 3 EV Aggregators} \label{CaseB}
We model our first coalition case study for 3 EV fleet aggregators.
Case B has the following properties:
\begin{itemize}
    \item 3 fleet aggregators with 3 different depots (the depot is the location of an EV before the first trip of the scheduling period);
    \item 5 charging stations, each fleet owns at least one station;
    \item Total number of customer nodes is 12 (these are temporary locations where an EV picks up or drops off a passenger).
\end{itemize}

The value function for each individual aggregator as well as the coalition value functions, together with energy consumption for each sub-coalition, are shown in Table~\ref{table:3 aggregator results}.
\begin{table}[ht]
\begin{center}
\caption{Case B: Coalitions Result Summary}
\begin{tabular}{r || c r r r | r}
    \hline
    Coalition & Value Fnc(\$) & $\varphi_1(\$)$ & $\varphi_2(\$)$ & $\varphi_3(\$)$ & E (kWh)   \\ 
    \hline\hline
        \{1\}     & 1.850 & 1.850     & -     & -     & 7.488   \\    \relax   
        \{2\}     & 3.044 & -         & 3.044 & -     & 12.321  \\    \relax  
        \{3\}     & 2.651 & -         & -     & 2.651 & 10.731  \\    \relax
        \{1, 2\}  & 1.635 & 0.220     & 1.414 & -     & 6.618   \\    \relax
        \{1, 3\}  & 1.599 & 0.399     & -     & 1.200 & 6.474   \\    \relax
        \{2, 3\}  & 2.187 & -         & 1.290 & 0.897 & 8.854   \\    \relax
        \{1, 2, 3\}& 1.381& {-}0.062  & 0.828 & 0.614 & 5.586   \\
        [1ex]
    \hline
\end{tabular}
\label{table:3 aggregator results}
\end{center}
\vspace{-0.5cm}
\end{table}
The cost of the grand coalition (last row on Table \ref{table:3 aggregator results})  are: $\varphi_1 = {-}0.062$, $\varphi_2 = 0.828$, and $\varphi_3 = 0.614$.
The aggregator 1 gets a negative cost which implies that part of the contribution of aggregators 2 and 3 will be given to 1. This is, however, not an aberration. Aggregator 1, working alone, has the least cost and thus will be profitable in any coalition which it joins. Thus, any other aggregator wants to partner with aggregator 1 because it will lead to a significant cost reduction.

\begin{figure}[ht]
    \centering
    \includegraphics[scale = 0.45]{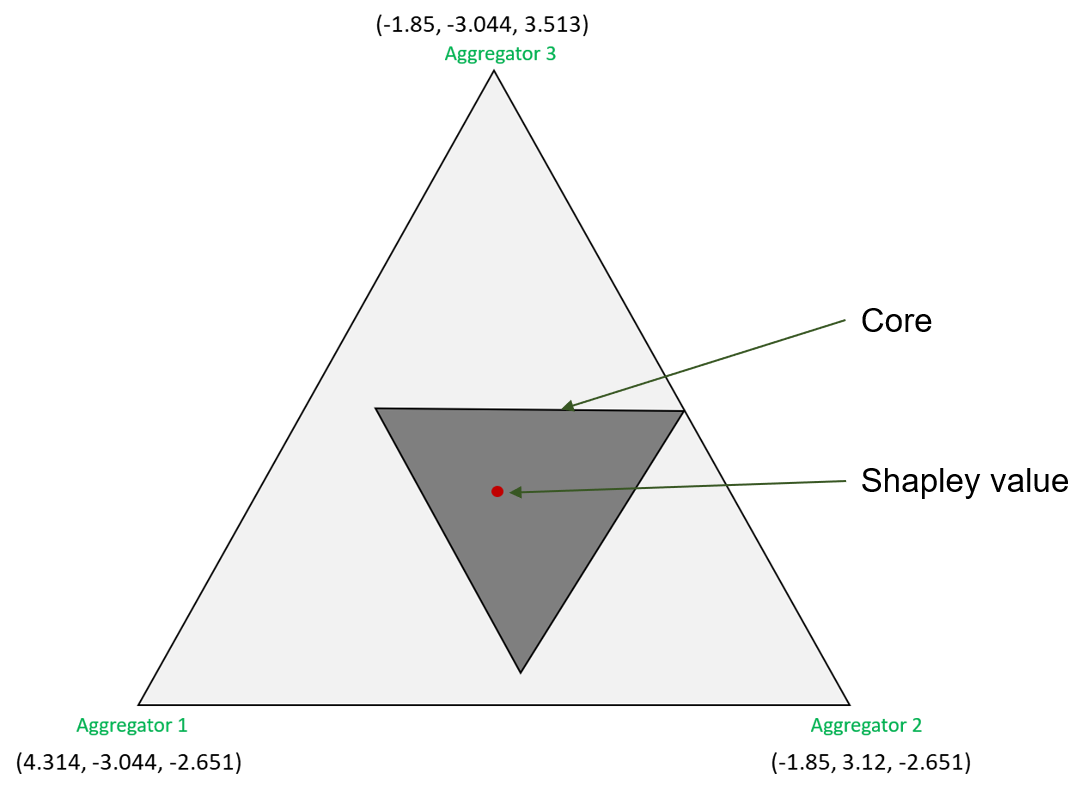}
    \caption{The Core and Shapley Value of a 3-aggregator Game}
    \label{fig:core3}
    \vspace{-0.5cm}
\end{figure}

From Table~\ref{table:3 aggregator results}, the best outcome for any aggregator is when all three cooperate. In the absence of this all-three coalition, aggregators 2 \& 3 have the incentive only to cooperate with each other, i.e. no other 2-aggregator coalition results in a Nash Equilibrium aside \{2, 3\}. The Core is shown in Fig~\ref{fig:core3} as the dark grey area while Shapley value is the red circle. The relative position of each player is shown in green text.

The Core shows all available coalitions with respective values that each aggregator might accept. Thus, we could have scenarios where concessions are made where aggregators do not accept the hard values given by Shapley computation. For example, what if aggregators $2$ and $3$ decide not to pay the negative cost to $1$ i.e. aggregator $1$ pays nothing? Aggregator $1$ would still be satisfied (see Fig~\ref{fig:core3} where the tip of aggregator 1 touches the axis). In fact, $1$ would be contented until his cost increases to \$0.4. This cost of \$0.4 for aggregator $1$ is the highest cost he can pay until $2$ \& $3$ have no additional incentive to refrain from the grand coalition.
As shown in the preceding paragraph, aggregators $2$ \& $3$ have an incentive to cut $1$ out of the whole cooperation framework and pursue their own coalition, \{2, 3\}. This incentive for $1$ to be left out could prompt it to accept costs set out by $2$ \& $3$. 
% All the acceptable costs by all players are shown in the dark gray area in Fig~\ref{fig:core3}. For the scenario where aggregator $1$ pays a higher cost than the Shapley value allocation, there is an increase in overall cost. This could either correspond to:
% \begin{enumerate}
%     \item A price of stability of 33\% over the original optimal grand coalition cost. The price of stability is the ratio of the new equilibrium cost to the minimum cost across all sub-coalitions. This 33\% which corresponds to the \$0.4 paid by aggregator 1 could be used to:
%     \begin{itemize}
%         \item set up an EV charging station co-managed by and accessible to all aggregators
%         \item create a low voltage feeder station exclusively for EV fleet belonging to all aggregators.
%     \end{itemize}  
    
%     \item A comparable 32\% decrease in cost paid by $2$ \& $3$
% \end{enumerate} 

Energy consumption follows the same trend as cost; single aggregators use more energy but forming coalitions leads to reduced energy consumption. Additionally, ability to use each others resources means vehicles do not have to travel long distances before their next recharge.
If we compare energy usage when aggregators are acting alone and the grand coalition, the maximum energy saving is observed between when \{2\} and \{1,2,3\}. This scenario leads to reduction in energy consumption by half while the least energy savings is a reduction of energy consumption of \{1\} by a quarter. This reduction in energy consumption directly translates to a reduction in emissions.
Here, we compare the two benchmark scenarios in Tables~\ref{table:3 aggs caseA1} \& \ref{table:3 aggs caseA2} with the solution from 3-Aggregators. If we compare case B with A1, the individual value functions and energy usages are somewhat similar. Like we remarked earlier, further interactions would lead to the aggregators settling on an equilibrium state. However, total energy savings when aggregators form the grand coalition is more significant. 

Case A2, however, mirrors the coalition results qualitatively, though in a rather counter-intuitive way. Table~\ref{table:3 aggregator results} shows that aggregator 2 working alone had the highest cost (higher cost means ``undesirable outcome'' in this case) and the greatest energy consumption, while Table~\ref{table:3 aggs caseA2} showed that Aggregator 2 had the lowest cost. In the aggregator-blind scenario, aggregators had a common belief space which was adequate knowledge of the transport and energy network. Thus, it was clear which aggregator was best situated for a particular routing decision. Aggregator 2 had the worst cost and participated in less routing decisions. The same argument goes for the grand coalition result shown in the last row of Table~\ref{table:3 aggregator results}. Actual routing and energy usage results showed that Aggregator 1 took the bulk share. Again, this is confirmed in the first row of Table~\ref{table:3 aggs caseA2}.

\subsection{Case C: 4 EV Aggregators} \label{CaseC}
Case C has the following properties:
\begin{itemize}
    \item Four fleet aggregators with 4 different depots; each fleet has the same number of vehicles as others. The fleets are also distributed strategically across the transport network such that aggregators 2, 3 \& 4 are centrally located while 1 is far from the center.
    \item Six charging stations, each fleet owning at least one;
    \item Total number of customer nodes is extended to 15 to afford more interaction between aggregators.
\end{itemize}

\begin{figure}[ht]
     \centering
     \includegraphics[scale = 0.5]{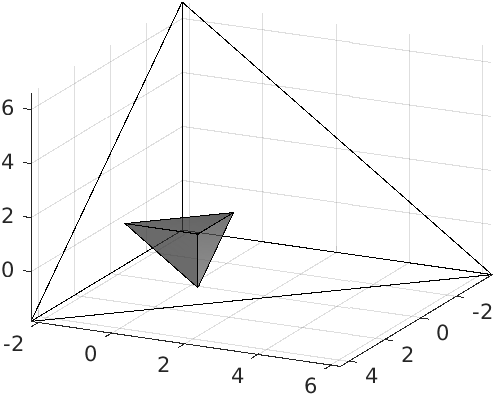}
     \caption{The Core of a 4-aggregator Game}
     \label{fig:core4}
\end{figure}

The Core is shown in Fig~\ref{fig:core4}. Unlike the case with 3 aggregators, the Core here is smaller and players do not have a big margin for deviation from the grand coalition. The Shapley value is given by the following vector:
$$[\varphi_1, \varphi_2, \varphi_3, \varphi_4] = [0.771, 0.176, 0.054, {-}0.075]$$

Again, the fourth aggregator gets a negative cost. This means part of the contributions from the other aggregators will go to the fourth one. This is so because aggregator 4 has the lowest cost when each of the aggregators acts alone. The Core is shown in Fig~\ref{fig:core4}.

% This cost is directly proportional to the energy consumption graph shown in Fig~\ref{fig:4aggregators result}. In fact, out of the 8 sub-coalitions in which aggregator 4 appears, a negative cost accrues half of the time.

% This negative cost of aggregator 4 only has a good utility in multi-player scenarios. The only cases where aggregator 4 contributes to significant cost reduction is in sub-coalitions \{1,2,4\}, \{1,3,4\} and the grand coalition \{1,2,3,4\}. In spite of this advantage to 4, there is a sub-coalition that has the tendency of deviating from the grand coalition; \{1,2,3\} because each aggregator, 1, 2 \& 3, does better here than anywhere else. Thus, if 4 is offered a cost increase, by other aggregators, up to \$0.5, it will be accepted. Aggregator 4 thus, has no incentive to deviate from the grand coalition. 

% As enumerated in Case B above, the price of stability in this case corresponds to a 62\% increase in cost over the grand coalition cost. This cost increase could be used either for infrastructural projects or subsidizing the costs for aggregators who have an incentive to deviate from grand coalition.

% In terms of energy consumption, the same trend of reduced consumption as seen in \ref{CaseB} is replicated here. Fig~\ref{fig:4aggregators result} shows the energy trend from single aggregators to the grand coalition. In the figure, greatest energy savings is observed when aggregator 1 is acting alone compared to the grand coalition. In that scenario, there is a 74\% reduction in energy usage whereas consumption is reduced by a third when aggregator 4 acting alone is compared with the grand coalition.

% \begin{figure}[ht]
%      \centering
%      \includegraphics[scale = 0.4]{figs/4aggregatorsenergy.PNG}
%      \caption{Energy Consumption per coalition of a 4-aggregator Game}
%      \label{fig:4aggregators result}
% \end{figure}

\subsection{Scalability Analysis}
\begin{figure}[ht]
     \centering
     \includegraphics[scale = 0.375]{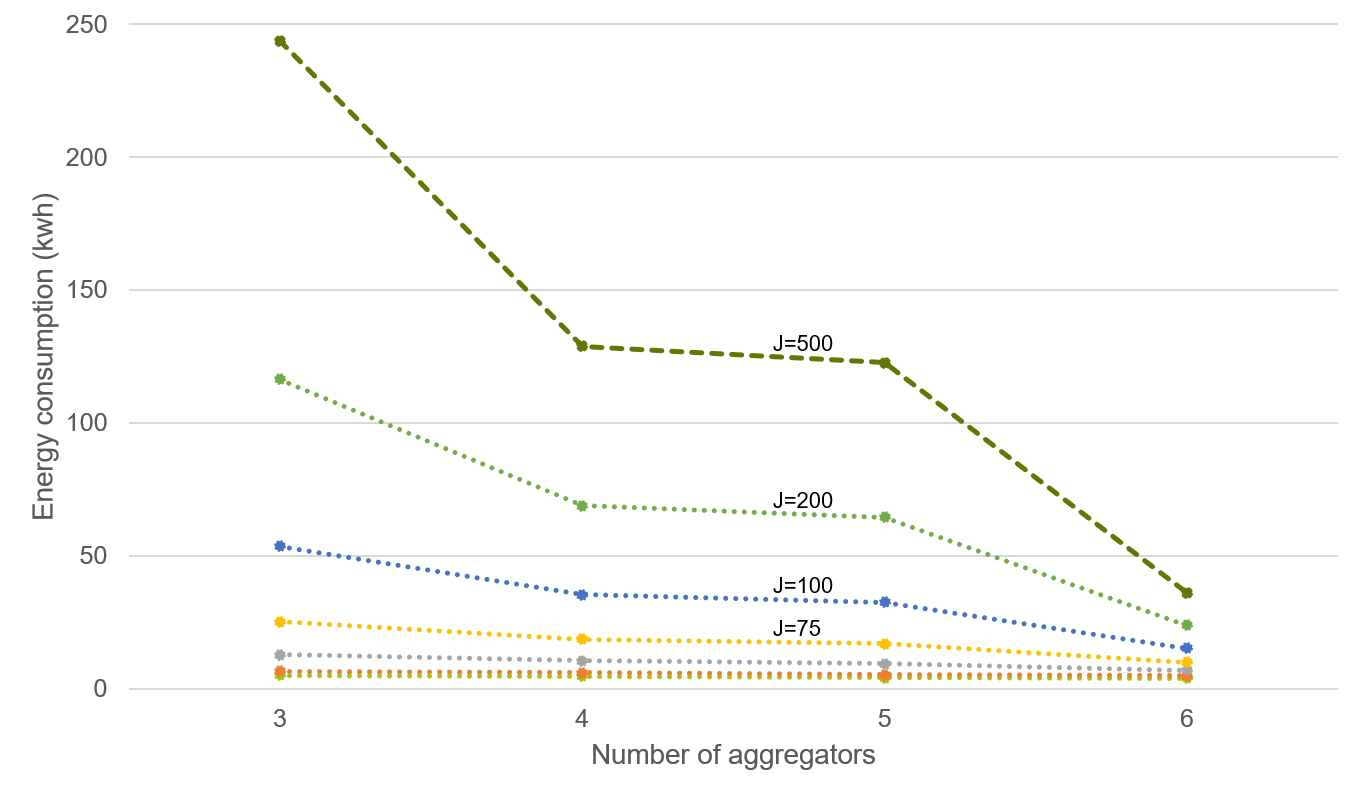}
     \caption{Energy consumption for increasing nodes}
     \vspace{-0.5cm}
     \label{fig:trends}
\end{figure}

Here, we carry out a sensitivity analysis by considering energy consumption trend with many nodes as well as many aggregators. As will be expected, many aggregators cooperating leads to reduced energy consumption.
In Fig~\ref{fig:trends}, the vertical axis represent the energy consumption over the timeline of routing and scheduling. The horizontal axis indicates the number of aggregators while $J$ represents the number of nodes in the Downtown Manhattan transport and energy network.

There is a general trend of decrease in cost and energy consumption as the number of aggregators rises. Thus, an increase of the economic efficiency. For a few locations on the network, this trend might mean a little decrease in cost or energy consumption. However, as the number of nodes increases, we observe a sharp decrease in energy usage and cost as aggregators increase. When there are multiple aggregators, there is a high probability that they would be spread out across the grid such that there can  be efficient service provision. This then leads to lesser travel time for individual vehicles.
In terms of energy requirement, the same argument can be made in this regard. As EVs are distributed uniformly across the grid, thereby reducing the distribution grid congestion.

\section{Conclusion} \label{conclusion}
This paper presents a cooperative game theory approach for evaluating the value of decentralized coordination of EV fleets by aggregators. A new EV-power grid integrated system model was introduced to describe the technical operation that our formulation was based on. We proposed two games with incomplete information, namely, aggregator-aware and aggregator-blind scenarios, as benchmarks for our coordination schemes.
Inter-fleet coordination was then addressed using cooperative game modelling. Allocation rules were set based on cost of operation of the grid and EVs while we used the Core to define the stability of the games. For efficiency and fairness, we utilized the Shapley value for cost allocation.

The objective of our cooperative game is the result of an EV transport routing problem that considers the distribution power grid. This objective also forms the marginal contribution of each aggregator to the coalitions/sub-coalitions which the aggregator is part of. We carried out tests on a set of nodes with 3 to 6 aggregators and were able to come to some conclusions:
\begin{enumerate}
    \item The relative location of a fleet aggregator influences his contribution to any coalition he is a part of; if a fleet owner is not close to pickup and drop-off of demands on the transport network, he would pay more than those fleet owners who are otherwise located.
    \item The coalition framework mirrors the aggregator blind scenario in terms of effectively reducing transportation cost and minimizing energy consumption by vehicle fleet.
    \item An aggregator can accrue negative cost; i.e. would receive extra payment from other aggregators. This might be the result of Shapley Value but the Core shows that other payment arrangements are possible.
    \item There is a general trend of decrease in cost and energy consumption as the number of aggregators rise. For few locations on the network, this trend might mean small percentage cost or energy consumption decrease. As the number of nodes increase, we observe a sharp decrease in energy usage and cost as aggregators increase. This could be a result of even distribution of aggregators across the network; leading to lesser travel time for individual vehicles and also less stringent demands on the distribution grid. 
\end{enumerate}

Future directions for this work may include investigating quantitatively ancillary services provision to the grid. The implication has been made here that extra cost accrued from cooperation might be utilized for this; however, how much realistically can these costs cater for. Another section not assessed in depth is the uncertainties in power consumption at local level. While power consumption at specific periods might be predictable (as adjudged from the daily load profiles of residences, commercial buildings and transport facilities), some level of stochasticity is still expected in power consumption. Thus, this could be a subject of further research.

\bibliographystyle{ieeetr}
\bibliography{refs}

\end{document}